\documentclass[11pt,leqno]{article}
\normalfont
\renewcommand{\baselinestretch}{1.25}

\usepackage{graphicx}
\usepackage{amsfonts}
\usepackage{mathrsfs}
\usepackage{rotating}
\usepackage{mathdots}
\usepackage{mathabx}
\usepackage[dvipsnames]{xcolor}

\newcommand{\IntroSection}{0}
\newcommand{\FibSection}{1}
\newcommand{\ContinuantSection}{2}
\newcommand{\BezoutSection}{3}
\newcommand{\EuclidSection}{4}
\newcommand{\ClosingSection}{5}

\newcommand{\FibTriangleFigure}{Figure 1.1}
\newcommand{\InterpretationOne}{Interpretation 1.1}
\newcommand{\InterpretationTwo}{Interpretation 1.2}
\newcommand{\Bijections}{Proposition 1.3}

\newcommand{\RecurrenceProp}{Proposition 2.1}

\newcommand{\BezoutLemma}{Lemma 3.1}

\newcommand{\EAExample}{Example 4.1}
\newcommand{\EATheorem}{Theorem 4.2}

\newfont{\mysmallscbolditalics}{ecoc0500 at 10pt}
\newcommand{\mymathfrak}[1]{\mbox{\mysmallscbolditalics #1}}

\newcommand{\QED}{\raisebox{0.5mm}{\fbox{\rule{0mm}{1.5mm}\ }}}

\newcounter{myfn}[page]
\renewcommand{\thefootnote}{\fnsymbol{footnote}}

\newcommand{\polyx}{\mbox{\large \sffamily x}}

\newcommand{\myt}{\mbox{\sffamily t}}

\begin{document}
\pagenumbering{arabic}
\thispagestyle{empty}%
\vspace*{-0.7in}
\hfill {\footnotesize March 10, 2023}

\begin{center}
{\large \bf Olry Terquem's forgotten problem} 

\vspace*{0.05in}
\renewcommand{\thefootnote}{1} 
Robert G.\ Donnelly,\footnote{Department of Mathematics and Statistics, Murray State University, Murray, KY 42071\\ 
\hspace*{0.25in}Email: {\tt rob.donnelly@murraystate.edu}}
\renewcommand{\thefootnote}{2}
\hspace*{-0.07in}Molly W.\ Dunkum,\footnote{Department of Mathematics, Western Kentucky University, Bowling Green, KY 42101\\ 
\hspace*{0.25in}Email: {\tt molly.dunkum@wku.edu}}
\renewcommand{\thefootnote}{3}
\hspace*{-0.07in}and Rachel McCoy\footnote{Department of Mathematics, Western Kentucky University, Bowling Green, KY 42101\\ 
\hspace*{0.25in}Email: {\tt mccoy.800@osu.edu}} 
\end{center} 

\vspace*{-0.17in}
\begin{abstract}
`Terquem's problem' is a name given in the twentieth century to the problem of enumerating certain integer sequences whose entries alternate in parity.  
In particular, this problem asks for the count of strictly increasing length $m$ sequences of positive integers bounded above by some integer $n$ whose odd-indexed entries are odd and whose even-indexed entries are even. 
This problem and its generalizations have been well-studied. 
However, the putative original source for this problem, an 1839 paper by Olry Terquem, is subtly different from the problem that is now attributed to Terquem. 
In this paper, we highlight this distinction and also make connections between Terquem's `forgotten' problem and the Fibonacci sequence, continuants, B\'{e}zout's Lemma, and the extended version of the Euclidean Algorithm.  

\begin{center}
\end{center} 
\end{abstract}

{\bf \S \IntroSection. Introduction.} 
Olry Terquem was a nineteenth century polymath -- a mathematician, historian, journal founder and editor, religious reformer, and officer of the French Legion of Honor. 
He is known for mathematical contributions in geometry, but in combinatorics his name is associated with an enumerative problem called `Terquem's problem'.   
Contemporary renderings of Terquem's problem generally follow the classical introductory combinatorics text by John Riordan, who seems to have borrowed it from an early twentieth century combinatorics text by the German mathematician Eugen Netto (see \S 49 in Ch.\ 3 of \cite{Netto}).  
Here is Riordan's version (\cite{Riordan}, Ch.\ 1 problem 15), with notation lightly modified to better fit our context:\\   
\hspace*{0.25in}\rule[-17mm]{0mm}{35.5mm}\parbox{5.95in}{\small Terquem's problem.  For combinations of $n$ numbered things in natural (rising) order, with $\myt(n,m)$ the number of $m$-combinations with odd elements in odd position and even elements in even positions, or, what is the same thing, with $\myt(n,m)$ the number of combinations with an equal number of odd and even elements for $m$ even and with the number of odd elements one greater than the number of even for $m$ odd, show that $\myt(n,m)$ has the recurrence $\myt(n,m) = \myt(n-1,m-1)+\myt(n-2,m)$ [with] $\myt(n,0)=1$, [that] $\myt(n,m)= {\lfloor \frac{n+m}{2} \rfloor \choose m}$, [and that] $f_{n+1} = \sum_{m=0}^{n}\myt(n,m) = f_{n} + f_{n-1}$.}\\ 
The objects to be enumerated in this version of Terquem's problem are called ``alternating parity subsets'' in the OEIS entry \cite{OEIS2}.    
(For example, when $n=4$, we get $\myt(4,0)=1$, $\myt(4,1)=2$, $\myt(4,2)=3$, $\myt(4,3)=1$, and $\myt(4,4)=1$, affording the partition $1+2+3+1+1$ of the Fibonacci number $8$.) 
Such alternating parity subsets have been well-studied and generalized; for a recent example, see the paper \cite{MM} by Toufik Mansour and Augustine O. Munagi and references therein. 

\renewcommand{\thefootnote}{4} 
However, the above version of Terquem's problem omits a requirement he imposed in his 1839 paper \cite{Terquem}.\footnote{During our research on this topic, the third-listed author helpfully translated Terquem's 1839 paper \cite{Terquem} (the original source is in French) as well as Netto's rendering of Terquem's problem in the 1901 text \cite{Netto} (the original source is in German).} 
That requirement was a {\em parity condition} on the lengths of the alternating parity subsets he considered. 
This somewhat subtle but forgotten difference in Terquem's original problem actually places his problem in close proximity to some very classical enumerative and number theoretic phenomena, namely, a right-triangular Fibonacci-type array of numbers, continuant polynomials, B\'{e}zout's Lemma, and the Euclidean extended algorithm. 
The main goal of this note is to revisit Terquem's original problem and to showcase these connections. 

\renewcommand{\thefootnote}{5} 
{\bf \S \FibSection. A Fibonacci-type array.} 
It is well-known that the sequence of sums of certain diagonals of Pascal's triangle is exactly the Fibonacci sequence $f_{0}=1$, $f_{1}=1$, $f_{2}=2$, $f_{3}=3$, $f_{4}=5$, etc.\footnote{This observation has been credited to Edouard Lucas, see Thomas Koshy's text \cite{Koshy}.} 
Indeed, we can skew Pascal's triangle as depicted below in order to view its Fibonacci diagonals as rows of the right-triangular array depicted in \FibTriangleFigure.   

\begin{figure}[ht]{
\begin{center}
\renewcommand{\baselinestretch}{1.2}
\footnotesize
\begin{tabular}{cccccccccc}
1 & & & & & & & & & \\
1 & & & & & & & & & \\
1 & 1 & & & & & & & & \\
1 & 2 & & & & & & & & \\
1 & 3 & 1 & & & & & & & \\
1 & 4 & 3 & & & & & & & \\
1 & 5 & 6 & 1 & & & & & & \\
1 & 6 & 10 & 4 & & & & & & \\
1 & 7 & 15 & 10 & 1 & & & & & \\
1 & 8 & 21 & 20 & 5 & & & & & \\
1 & 9 & 28 & 35 & 15 & 1 & & & & \\
1 & 10 & 36 & 56 & 35 & 6 & & & & \\
1 & 11 & 45 & 84 & 70 & 21 & 1 & & & \\
1 & 12 & 55 & 120 & 126 & 56 & 7 & & & \\
1 & 13 & 66 & 165 & 210 & 126 & 28 & 1 & & \\
1 & 14 & 78 & 220 & 330 & 252 & 84 & 8 & & \\
1 & 15 & 91 & 286 & 495 & 462 & 210 & 36 & 1 & \\
1 & 16 & 105 & 364 & 715 & 792 & 462 & 120 & 9 & 
\end{tabular}

ETC.

{\bf \FibTriangleFigure}
\vspace*{-0.25in}
\end{center}}
\end{figure}

\renewcommand{\thefootnote}{6} 
\noindent 
We call this the {\em Fibonacci right-triangular array}, or simply the {\em Fibonacci array}. 
We refer to this triangular array using the notation $\mymathfrak{F}$ and label the $k^{\mbox{\tiny \em th}}$ entry of the $n^{\mbox{\tiny \em th}}$ row as $f_{n,k}$, so that \[\mymathfrak{F} = (f_{n,k})_{n\in\{0,1,2,\ldots\}, k\in\{0,1,\ldots,\lfloor\frac{n}{2}\rfloor\}}.\] 
Formally, we build $\mymathfrak{F}$ by declaring $f_{0,0} := f_{1,0} := 1$ and then applying the recurrence 
\[f_{n,k} := f_{n-1,k} + f_{n-2,k-1}\] 
for integers $k$ and $n$ with $n \geq 2$ and with the understanding that $f_{n,k} := 0$ when $k < 0$ or $k > \lfloor \frac{n}{2} \rfloor$. 
Using this defining recurrence, one can easily check that for any integers $n$ and $k$ with $0 \leq k \leq \lfloor \frac{n}{2} \rfloor$ we have
\[f_{n,k} = {n-k \choose k} \hspace*{0.25in}\mbox{and}\hspace*{0.25in}\sum_{i=0}^{\lfloor \frac{n}{2} \rfloor}f_{n,i} = f_{n}.\] 
So, the sequence of entries on the $n^{\mbox{\tiny \em th}}$ row of $\mymathfrak{F}$ can be thought of as a refinement of the $n^{\mbox{\tiny \em th}}$ Fibonacci number. 
Indeed, some collections of polynomials, such as those studied by Lucas\footnote{On p.\ 186 of \cite{Lucas}, the quantity $U_{n}$ (where $n$ is a positive integer) is a polynomial in the integer variables $P$ and $Q$ and has sign-alternating coefficients from row $n-1$ of $\mymathfrak{F}$.} (see \cite{Lucas}) and Ernst Jacobsthal (see \cite{Koshy} Ch.\ 39), have coefficients from the rows of $\mymathfrak{F}$. 
The numbers comprising these refinements also have many enumerative interpretations. 
Next, we present two such interpretations -- one well-known and one, due to Olry Terquem, that seems to have been forgotten. 

\noindent
{\bf \InterpretationOne: A classical enumerative interpretation of the Fibonacci array.}  Fix a positive integer $n$ and a nonnegative integer $m$.  
A {\em consecutive-free} sequence $S = (S_{1}, S_{2}, \ldots, S_{m})$ {\em strictly bounded by} $n$ is a strictly increasing sequence of $m$ positive integers wherein $S_{m} < n$ and $S_{i+1} - S_{i} \geq 2$ for any $i \in \{1,2,\ldots,m-1\}$. 
We say that $S$ has length $m$. 
When $m=0$, then $S$ is the empty sequence, is denoted by $()$, is considered to have length zero, and is considered to be a consecutive-free sequence. 
Let $\mathcal{CF}(n)$ be the set of consecutive-free sequences strictly bounded by $n$. 
So, for example, 
\begin{eqnarray*}
\mathcal{CF}(5) & = & \{(), (1), (2), (3), (4), (13), (14), (24)\}\\
\mathcal{CF}(6) & = & \{(), (1), (2), (3), (4), (5), (13), (14), (15), (24), (25), (35), (135)\} 
\end{eqnarray*}
Further, let $\mathcal{CF}(n,k)$ denote the set of sequences in $\mathcal{CF}(n)$ having length $k$, when $k$ is a nonnegative integer. 
According to OEIS entry A011973 \cite{OEIS}, we have $|\mathcal{CF}(n,k)| = f_{n,k}$ for all positive integers $n$ and all integers $k$ with $0 \leq k \leq \lfloor \frac{n}{2} \rfloor$. 
In particular, $|\mathcal{CF}(n)| = f_{n}$, the $n^{\mbox{\em \tiny th}}$ Fibonacci number.\hfill\QED  

\noindent
{\bf \InterpretationTwo: Olry Terquem's forgotten problem.}  Again, assume $n$ is a positive integer and that $m$ is a nonnegative integer. 
An {\em oe-sequence} $T := (T_{1}, T_{2}, \ldots, T_{m})$ {\em bounded by} $n$ is a strictly increasing sequence of $m$ positive integers wherein $m$ and $n$ have the same parity, $T_{1}$ is odd, $T_{m} \leq n$, and $T_{i}$ and $T_{i+1}$ have opposite parity for any $i \in \{1,2,\ldots,m-1\}$. 
(In this definition, the ``oe'' prefix is meant to suggest the alternating odd-even pattern of an oe-sequence.) 
We say that $T$ has length $m$ and write $\ell(T) = m$. 
When $m=0$, then $T$ is the empty sequence, is denoted by $()$, is considered to have length zero, and is considered to be an oe-sequence. 
Let $\mathcal{OE}(n)$ be the set of oe-sequences bounded by $n$. 
For reasons to be made apparent shortly, we set $\mathcal{OE}(-1) := \emptyset$ and $\mathcal{OE}(0) := \{()\}$. 
So, for example,
\begin{eqnarray*}
\mathcal{OE}(5) & = & \{(1), (3), (5), (123), (125), (145), (345), (12345)\}\\
\mathcal{OE}(6) & = & \{(), (12), (14), (16), (34), (36), (56), (1234), (1236), (1256), (1456), (3456), (123456)\} 
\end{eqnarray*}
Further, for $0 \leq k \leq \lfloor \frac{n}{2} \rfloor$, let $\mathcal{OE}(n,k)$ be the set of sequences $T$ in $\mathcal{OE}(n)$ with $\ell(T) = n - 2k$. 
{\em Terquem's forgotten problem} is to show that $|\mathcal{OE}(n,k)| = f_{n,k}$. 
In \Bijections\ below, we solve this problem by demonstrating a bijection between $\mathcal{CF}(n,k)$ and $\mathcal{OE}(n,k)$. 

Likewise, an {\em eo-sequence} has all the properties of an oe-sequence except that it begins with an even integer instead of an odd and its length has parity opposite that of $n$;  
we then define $\mathcal{EO}(n)$ to be the set of eo-sequences bounded by $n$. 
This time we set $\mathcal{EO}(-1) := \{()\}$ and $\mathcal{EO}(0) := \emptyset$. 
We have, for example,
\begin{eqnarray*}
\mathcal{EO}(5) & = & \{(),(23),(25),(45),(2345)\}\\
\mathcal{EO}(6) & = & \{(2), (4), (6), (234), (236), (256), (456), (23456)\} 
\end{eqnarray*}
From here on, we refer to oe- and eo-sequences generically as {\em oereo sequences} (pronounced `OH-ree-oh'). 
For $0 \leq k \leq \lfloor \frac{n-1}{2} \rfloor$, let $\mathcal{EO}(n,k)$ be the set of $T \in \mathcal{EO}(n)$ with $\ell(T) = (n-1) - 2k$. 
In \Bijections\ below, we demonstrate a bijection between $\mathcal{EO}(n,k)$ and $\mathcal{OE}(n-1,k)$ and conclude that $|\mathcal{EO}(n,k)| = f_{n-1,k}$. 
\hfill\QED

To set up our solution to Terquem's forgotten problem in \Bijections, we require the following ideas. 
When entries of a given sequence are distinct, we may regard the sequence as a set and apply set operations. 
So, when $n$ is a positive integer, we observe that a strictly increasing sequence $(T_{1},\ldots,T_{m})$ of integers with $1 \leq T_{1}$ and $T_{m} \leq n$ is an oe-sequence if and only if $\{1,2,\ldots,n\} \setminus \{T_{1},\ldots,T_{m}\}$ is the union of even-length subsequences of consecutive integers. 

For a nonnegative integer $k \leq \lfloor \frac{n}{2} \rfloor$, we define a function $\phi: \mathcal{CF}(n,k) \longrightarrow \mathcal{OE}(n,k)$ in the following way.  
Given $S = (S_{1},\ldots,S_{k}) \in \mathcal{CF}(n,k)$, let 
\[T := \phi(S) := \left\{\rule[-1.75mm]{0mm}{4.5mm}1,2,\ldots,n-1,n\right\} \mbox{\LARGE $\setminus$} \left(\mbox{\footnotesize $\bigcup\limits_{i=1}^{k}$} \{S_{i},S_{i}+1\}\right).\] 
Clearly $\ell(T) = n-2k$, and therefore $\ell(T)$ has the same parity as $n$.  
Since $T$ is formed from $\{1,2,\ldots,n\}$ by removing consecutive pairs of integers, then the set of numbers between any two successive entries of $T$ is an even-length (and possibly length zero) string of consecutive integers; so, as observed in the previous paragraph, $T$ is an oe-sequence. 
Next, we define a function $\psi: \mathcal{EO}(n,k) \longrightarrow \mathcal{OE}(n-1,k)$ as follows. 
Declare that $() \stackrel{\psi}{\longmapsto} ()$.  
Otherwise, when $T=(T_{1},\ldots,T_{m})$ with $m = \ell(T) > 0$, we set $\psi(T) := (T_{1}-1,T_{2}-1,\ldots,T_{m}-1)$. 
It is evident that $\psi(T) \in \mathcal{OE}(n-1,k)$. 

\noindent 
{\bf \Bijections}\ \ {\sl Let $n$ and $k$ be integers with $n$ positive and $0 \leq k \leq \lfloor \frac{n}{2} \rfloor$.  Then the functions $\phi: \mathcal{CF}(n,k) \longrightarrow \mathcal{OE}(n,k)$ and $\psi: \mathcal{EO}(n,k) \longrightarrow \mathcal{OE}(n-1,k)$ are bijections.
Therefore $|\mathcal{OE}(n,k)| = f_{n,k}$ and $|\mathcal{EO}(n,k)| = f_{n-1,k}$.} 

{\em Proof.} It is observed above that the range of each of $\phi$ and $\psi$ is within the declared target set. 
Since the procedures employed by each of $\phi$ and $\psi$ are clearly reversible, we conclude that these functions are bijections. 
The remaining claims follow from the fact that $|\mathcal{CF}(n,k)| = f_{n,k}$, cf.\ \InterpretationOne\ above.\hfill\QED

{\bf \S \ContinuantSection. A combinatorial characterization of continuant polynomials.} 
We now define certain multivariate polynomials $g_{n}$ and $h_{n}$ as `generating functions' over certain sets of oereo sequences. 
These definitions are effected by the following general notion: For any strictly increasing positive integer sequence $T = (T_{1},\ldots,T_{m})$, set $\polyx_{T} := x_{T_{1}}x_{T_{2}} \cdots x_{T_{m}}$ (where each $x_{i}$ is an indeterminate), with $\polyx_{()} := 1$ for the empty sequence. 
For any integer $n\geq -1$, set
\[g_{n} := \sum_{T \in \mathcal{EO}(n)}\polyx_{T} \hspace*{0.3in}\mbox{and}\hspace*{0.3in} h_{n} :=  \sum_{T \in \mathcal{OE}(n)}\polyx_{T},\]
where of course an empty sum is zero. 
Note that each of $g_{n}$ and $h_{n}$ is a polynomial in $x_{1},\ldots,x_{n}$ with unit coefficients. 
Call $g_{n}$ an {\em eo-polynomial} and $h_{n}$ an {\em oe-polynomial}; together, these are generically referred to as {\em oereo polynomials}. 
For the record, here are the first several eo- and oe-polynomials:

\vspace*{-0.1in}
\begin{center}
{\small \begin{tabular}{rclcrcl}
$g_{-1}$ & = & $1$ & \hspace*{0.1in} & $h_{-1}$ & = & $0$\\
$g_{0}$ & = & $0$ & \hspace*{0.1in} & $h_{0}$ & = & $1$\\
$g_{1}$ & = & $1$ & \hspace*{0.1in} & $h_{1}$ & = & $x_{1}$\\
$g_{2}$ & = & $x_{2}$ & \hspace*{0.1in} & $h_{2}$ & = & $1+x_{1}x_{2}$\\
$g_{3}$ & = & $1+x_{2}x_{3}$ & \hspace*{0.1in} & $h_{3}$ & = & $x_{1}+x_{3}+x_{1}x_{2}x_{3}$\\
$g_{4}$ & = & $x_{2}+x_{4}+x_{2}x_{3}x_{4}$ & \hspace*{0.1in} & $h_{4}$ & = & $1+x_{1}x_{2}+x_{1}x_{4}+x_{3}x_{4}+x_{1}x_{2}x_{3}x_{4}$\\
$g_{5}$ & = & $1+x_{2}x_{3}+x_{2}x_{5}+x_{4}x_{5}+x_{2}x_{3}x_{4}x_{5}$ & \hspace*{0.1in} & $h_{5}$ & = & $x_{1}+x_{3}+x_{5}+ x_{1}x_{2}x_{3}+x_{1}x_{2}x_{5}$\\
 & & & \hspace*{0.15in} & & & \hspace*{0.25in}$+x_{1}x_{4}x_{5}+x_{3}x_{4}x_{5}+x_{1}x_{2}x_{3}x_{4}x_{5}$
\end{tabular}}
\end{center}

A consequence of \Bijections\  is that when $n \geq 1$, the 
Fibonacci number $f_{n-1}$ counts the number of terms in the polynomial $g_{n}$, and the 
Fibonacci number $f_{n}$ counts the number of terms in $h_{n}$. 
Moreover, the number of terms of degree $(n-1)-2k$ (respectively $n-2k$) for $0 \leq k \leq \lfloor \frac{n-1}{2} \rfloor$ (resp.\ $0 \leq k \leq \lfloor \frac{n}{2} \rfloor$) within $g_{n}$ (resp.\ $h_{n}$) is $f_{n-1,k}$ (resp.\ $f_{n,k}$). 
Given that these polynomials have an underlying Fibonacci-like structure, it is not surprising that they can be constructed using simple recurrence relations, as in our next result.  
We use the notation $g_{n}(c_{1},\ldots,c_{n})$ to indicate that each $x_{j}$ has been specialized to some given integer quantity $c_{j}$, and similarly for $h_{n}$.  

\noindent
{\bf \RecurrenceProp}\ \ {\sl The eo-polynomials $\{g_{n}\}_{n \geq -1}$ are uniquely determined by the recurrence relations $g_{-1} = 1$, $g_{0} = 0$, and $g_{n+1} = g_{n-1} + g_{n}x_{n+1}$ for $n \geq 0$.  
Similarly, the oe-polynomials $\{h_{n}\}_{n \geq -1}$ are uniquely determined by the recurrence relations $h_{-1} = 0$, $h_{0} = 1$, and $h_{n+1} = h_{n-1} + h_{n}x_{n+1}$ for $n \geq 0$. 
For a positive integer $n$ and real numbers $c_{1},\ldots,c_{n}$, we have $h_{n-1}(c_{2},\ldots,c_{n}) = g_{n}(c_{1},c_{2},\ldots,c_{n})$ and $c_{1}h_{n-1}(c_{2},\ldots,c_{n})+g_{n-1}(c_{2},\ldots,c_{n}) = h_{n}(c_{1},c_{2},\ldots,c_{n})$.} 

{\em Proof.} The claims in the first two sentences of the proposition follow from routine induction arguments, which we leave to the reader. 
The claims of the last sentence follow from the definitions of $g_{n}$ and $h_{n}$ in terms of eo- and oe-sequences. \hfill\QED

The recurrence in the preceding proposition allows us to identify oe-polynomials as the `continuant polynomials' of \S 4.5.3 of \cite{Knuth}. 
There, Knuth shows how such polynomials are related to continued fractions (hence their name) and to the Euclidean Algorithm. 
Knuth attributes to Euler (presumably referring to \cite{Euler}) the following characterization of the terms of the oe-polynomial $h_{n}$: They are obtained from $x_{1}x_{2}\cdots x_{n}$ by ``deleting zero or more nonoverlapping pairs of consecutive variables $x_{j}x_{j+1}$.''  
As observed in \S \FibSection\ above, these are exactly the oe-sequences bounded by $n$.

{\bf \S \BezoutSection. B\'{e}zout's Lemma -- a nonconstructive result?} 
The result traditionally known as B\'{e}zout's Lemma is the following assertion: 

\noindent 
{\bf \BezoutLemma}\ \  {\sl Given integers $a$ and $b$, not both of which are zero, there exist integers $s$ and $t$ such that $\gcd(a,b) = as+bt$.} 

Although this result is attributed to B\'{e}zout \cite{Bezout}, a version of it was discovered more than a century prior by Claude Gaspard Bachet \cite{Bachet}.  
Working within the context of the Euclidean Algorithm from Euclid's Book VII, Bachet used symbolic algebra to describe a general solution to the Diophantine equation $ax-by=1$, when $a$ and $b$ are co-prime (i.e.\ $\gcd(a,b)=1$), by backtracking through the quotients and remainders produced by the steps of the algorithm. 
But it was with B\'{e}zout's work that the more general notion of solving $ax+by = \gcd(a,b)$ over the integers was recognized as a key result within elementary number theory (see, for example, Maarten Bullynck's wonderfully readable account of the European history of such remainder problems in the centuries preceding Gauss \cite{Bullynck}). 

But, what might be practical about B\'{e}zout's Lemma? 
Here is a modern application. 
Cryptosystems often utilize modular arithmetic and, in particular, can require the computation of multiplicative inverses within this environment. 
Suppose, for example, we know that $a$ and $b$ are co-prime and wish to compute $\overline{b}^{-1}$ within the quotient ring $\mathbb{Z}_{a} := \mathbb{Z}/a\mathbb{Z}$, if such an inverse exists at all. 
(Here, $\overline{x}$ denotes the coset $x+a\mathbb{Z}$ in $\mathbb{Z}_{a}$.) 
Well, write $1 = as+bt$, where $(s,t)$ is an integer pair of `B\'{e}zout coefficients'.   
Then $\overline{1} = \overline{as+bt} = \overline{a}\, \overline{s} + \overline{b}\, \overline{t} = \overline{b}\, \overline{t}$, so $\overline{b}\, \overline{t} = \overline{1}$. 
That is, $\overline{b}$ has a multiplicative inverse in $\mathbb{Z}_{a}$, and $\overline{b}^{-1} = \overline{t}$. 

\renewcommand{\thefootnote}{7} 
Given such computational contexts, it is curious that contemporary treatments of B\'{e}zout's Lemma are sometimes nonconstructive and precede (as in \S 3.3 of \cite{Rosen}), or even obviate (as in Ch.\ 0 of \cite{Gallian}), the Euclidean Algorithm altogether. 
In such approaches, B\'{e}zout's Lemma is proved using the (nonconstructive) Well Ordering Principle of the positive integers, which asserts that {\sl any nonempty set of positive integers contains a smallest integer.}\footnote{It is worth noting that in these same references, the Division Algorithm is also taken to be a nonconstructive result whose existence aspect is a consequence of the Well Ordering Principle (see Ch.\ 0 of \cite{Gallian} and \S 1.5 of \cite{Rosen}). In particular, when $a$ and $b$ are integers with $b \ne 0$, then $\{a-bk\, |\, k \in \mathbb{Z} \mbox{ and } a-bk \geq 0\}$ is nonempty and contains a smallest element $r$ of the form $a-bq$, from which we easily deduce that $0 \leq r < |b|$. Of course, such nonconstructive reasoning is not at all necessary in order to settle the existence claims of the Division Algorithm.}
The proof goes as follows: Assuming $a$ and $b$ are not both zero, then $\{ax+by\, |\, x,y \in \mathbb{Z} \mbox{ and } ax+by>0\}$ is nonempty and we may, by invoking the Well Ordering Principle, declare $d$ to be the smallest element in this set.  
Clearly $\gcd(a,b)|d$, so $\gcd(a,b) \leq d$. 
To see that $d$ is a divisor of $a$, write $a = dq+r$ with $0 \leq r < d$, whence $r = dq-a = a(xq-1)+b(yq)$, forcing $r=0$. 
We similarly see that $d|b$, and we may conclude now that $d \leq \gcd(a,b)$. 
So, $d=\gcd(a,b)$. 

Shortly, we will show how oereo polynomials, together with data provided by the Euclidean Algorithm, afford explicit and non-recursive formulas for B\'{e}zout coefficients. 

\renewcommand{\thefootnote}{8} 
{\bf \S \EuclidSection. A version of the extended Euclidean Algorithm.} 
The Euclidean Algorithm is an ancient procedure for computing the greatest common divisor of two integers without factoring them.  
This algorithm is Euclid's in that it famously appears in his great treatise, the {\em Elements}, as Proposition 2 of Book VII. 
However, most likely this result was known before Euclid, perhaps by the Pythagoreans, as contended by (among others) the mathematician and historian Bartel Leendert van der Waerden\footnote{From p.\ 115 of \cite{vdW}: 
``... Book VII is not a later reconstruction, but a piece of ancient mathematics. ... Thus we have acquired an important insight: {\em Book VII was a textbook on the elements of the Theory of Numbers, in use in the Pythagorean school.}''} in \cite{vdW}. 
In his canonical work on {\em The Art of Computer Programming}, Donald Knuth writes that the Euclidean Algorithm is ``the granddaddy of all algorithms, because it is the oldest nontrivial algorithm that has survived to the present day'' (see \S 4.5.2 of \cite{Knuth}).  

\renewcommand{\thefootnote}{9} 
An extended version of the Euclidean Algorithm (e.g.\ \S 3.4 of \cite{Rosen}) determines B\'{e}zout coefficients along with the greatest common divisor.  
Here, we will make connections with B\'{e}zout's Lemma and our solution to Terquem's forgotten problem by formulating the (traditional, unextended) Euclidean Algorithm as a pseudo-code procedure that advances in steps. 
We take as our input a pair of positive\footnote{Positivity is assumed for the sake of simplicity.} integers $a$ and $b$ with $a \geq b$.  
The procedure applies the Division Algorithm in successive steps, producing at each step a quotient and a remainder, and returns the last nonzero remainder as the greatest common divisor. 
We will also return other data provided by the procedure, in particular, the complete list quotients and complete list of remainders obtained along the way. 
Using oereo polynomials, we will see how to recover from this data any information that is yielded by the extended Euclidean Algorithm. 

\noindent 
\centerline{\underline{\sc Euclidean Algorithm (EA)}}\\
\underline{\sc Input:} Positive integers $a$ and $b$ with $a \geq b$.\\
\underline{\sc Initialization:} Set $r_{-1} := a$ and $r_{0} := b$, to be regarded as {\em de facto} remainders. Let $i$ be our step-counting index, initialized with $i := 1$.\\
\underline{\sc Step 1:} Applying the Division Algorithm, we obtain unique integers $q_{i}$ and $r_{i}$ wherein $0 \leq r_{i} < r_{i-1}$ and 
\[r_{i-2} = r_{i-1}q_{i}+r_{i}.\]
\underline{\sc Step 2:} If $r_{i} = 0$, then {\sc Return} the following data: $\mbox{\sc gcd} := r_{i-1}$, $\mbox{\sc num{\_}steps} := i$, $\mbox{\sc rem{\_}list} := (r_{-1}, r_{0}, \ldots, r_{i})$, and $\mbox{\sc quo{\_}list} := (q_{1}, \ldots, q_{i})$.\\ 
\underline{\sc Step 3:} Otherwise, with $r_{i} > 0$, replace $i$ with $i+1$ and return to {\sc Step 1}.\\
\underline{\sc Output:} A positive integer $\mbox{\sc gcd}$ (that will be shown to be the greatest common divisor of $a$ and $b$), a positive integer $\mbox{\sc num{\_}steps}$ that records the number of times we apply the Division Algorithm, the list $\mbox{\sc rem{\_}list}$ of nonnegative integer remainders indexed from $-1$ up to $\mbox{\sc num{\_}steps}$, and the list 
$\mbox{\sc quo{\_}list}$ of positive integer quotients indexed from $1$ through $\mbox{\sc num{\_}steps}$.\\
\centerline{$\overline{\hspace*{4in}}$}

For positive integers $a$ and $b$ with $a \geq b$, we use the notation $\mbox{\sc EA($a,b$)[gcd]}$ to denote the greatest common divisor that the algorithm returns, $\mbox{\sc EA($a,b$)[num{\_}steps]}$ to denote the number of times EA applies the Division Algorithm, $\mbox{\sc EA($a,b$)[rem{\_}list]}$ to denote the list of remainders, and $\mbox{\sc EA($a,b$)[quo{\_}list]}$ to denote the list of quotients. 

\noindent 
{\bf \EAExample}\ \ 
For purposes of illustrating this version of the Euclidean Algorithm, we sought input numbers $a$ and $b$ such that the algorithm would terminate in exactly 7 steps, i.e.\ 7 applications of the Division Algorithm. 
How we found such $a$ and $b$ -- and why we wanted this exact number of steps -- will be revealed shortly. 
With this in mind, consider $a = 4449$ and $b=935$.  Here's how $\mbox{\sc EA}(a,b)$ works:
\begin{center}
\begin{tabular}{rclcrcl}
$4449$ & $=$ & $935 \cdot 4 + 709$ & \hspace*{1in} & $31$ & $=$ & $9 \cdot 3 + 4$\\
$935$ & $=$ & $709 \cdot 1 + 226$ & \hspace*{1in} & $9$ & $=$ & $4 \cdot 2 + 1$\\
$709$ & $=$ & $226 \cdot 3 + 31$ & \hspace*{1in} & $4$ & $=$ & $1 \cdot 4 + 0$\\
$226$ & $=$ & $31 \cdot 7 + 9$ & \hspace*{1in} & & & \\
\end{tabular}
\end{center}
EA output is as follows: $\mbox{\sc EA($a,b$)[gcd]} = 1$, $\mbox{\sc EA($a,b$)[num{\_}steps]} = 7$, $\mbox{\sc EA($a,b$)[rem{\_}list]} = (4449,935,709,226,31,9,4,1,0)$, and $\mbox{\sc EA($a,b$)[quo{\_}list]} = (4,1,3,7,3,2,4)$.  
In fact, we can construct B\'{e}zout coefficients by backtracking through the preceding equations: 
$1 = 9+4\cdot^{\mbox{\tiny $-$}}\!{2} = 31\cdot^{\mbox{\tiny $-$}}\!{2} + 9\cdot{7} = 226\cdot{7}+31\cdot^{\mbox{\tiny $-$}}\!{51} = 709\cdot^{\mbox{\tiny $-$}}\!{51}+226\cdot{160} = 935\cdot{160}+709\cdot^{\mbox{\tiny $-$}}\!{211} = 4449\cdot^{\mbox{\tiny $-$}}\!{211}+935\cdot{1004}$, so our B\'{e}zout coefficients are $s = ^{\mbox{\tiny $-$}}\!{211}$ and $t=1004$.\hfill\QED

\noindent
{\bf \EATheorem\ (A Theorem on the Euclidean Algorithm)}\ \ {\sl Let $a$ and $b$ be positive integers with $a \geq b$.}

{\sl (0) } \parbox[t]{6.05in}{{\sl With this input, the Euclidean Algorithm terminates. For the remainder of the theorem statement, we set}\\ 
\centerline{$G := \mbox{\sc EA($a,b$)[gcd]}$, $n := \mbox{\sc EA($a,b$)[num{\_}steps]}$,} 
\centerline{$(r_{-1}, r_{0}, \ldots r_{n}) := \mbox{\sc EA($a,b$)[rem{\_}list]}$, {\sl and} $(q_{1},\ldots,q_{n}) := \mbox{\sc EA($a,b$)[quo{\_}list]}$,}}

\vspace*{0.05in}
\hspace*{0.25in} \parbox[t]{6.05in}{{\sl where (as claimed at the conclusion of the algorithm) the output $G=r_{n-1}$ is a positive integer, $r_{-1}$ and $r_{0}$ are the positive integers $a$ and $b$ respectively, $n$ records the number of times the Euclidean Algorithm applies the Division Algorithm, and $(r_{1}, \ldots r_{n})$ is the 
sequence of nonnegative integer remainders and $(q_{1},\ldots,q_{n})$ the sequence of positive integer quotients obtained from these $n$ applications of the Division Algorithm.}}

\vspace*{0.1in}
{\sl (1) } \parbox[t]{6.05in}{{\sl (a B\'{e}zout-type result) For any $i \in \{-1,0,\ldots,n\}$, set $\overrightarrow{g_{i}} := g_{i}(q_{1},\ldots,q_{i})$, $\overleftarrow{g_{i}} := g_{i}(q_{n+1-i},\ldots,q_{n})$, $\overrightarrow{h_{i}} := h_{i}(q_{1},\ldots,q_{i})$,  and $\overleftarrow{h_{i}} := h_{i}(q_{n+1-i},\ldots,q_{n})$, and momentarily allow $r_{n+1} := G$. For each such $i$, $r_{i} = (-1)^{i+1}\overrightarrow{g_{i}}\, a + (-1)^{i}\overrightarrow{h_{i}}\, b${\huge ,}\,  $r_{n-i} = \overleftarrow{g_{i}}\, G${\huge ,}\,  and $r_{n-1-i} = \overleftarrow{h_{i}}\, G$.}} 

{\sl (2) } \parbox[t]{6.05in}{{\sl For any $i \in \{0,1,\ldots,n\}$, if $c$ is a divisor of $r_{i-1}$ and $r_{i}$, then $c$ divides each of $r_{-1}$, $r_{0}$, $\ldots$, $r_{n}$.}} 

{\sl (3) } \parbox[t]{6.05in}{{\sl The greatest common divisor of $a$ and $b$ is the quantity $G$, i.e.\ $G = \gcd(a,b)$. 
Moreover, the set $\{ax+by\, |\, x,y \in \mathbb{Z}\}$ consists exactly of the integer multiples of $\gcd(a,b)$. 
In particular, $\gcd(a,b)$ is the smallest positive integer expressible as an integer linear combination of $a$ and $b$. 
Moreover, $\gcd(a,b) = (-1)^{n}g_{n-1}(q_{1},\ldots,q_{n-1})a+(-1)^{n-1}h_{n-1}(q_{1},\ldots,q_{n-1})b${\huge ,}\, $b = g_{n}(q_{1},\ldots,q_{n}){\gcd(a,b)}${\huge ,}\, and $a = h_{n}(q_{1},\ldots,q_{n}){\gcd(a,b)}$.}}

{\em Proof}. For {\sl (0)}, consider the sequence of remainders $r_{-1} \geq r_{0} > r_{1} > \cdots$.   
Now, $b = r_{0} > 0$, $b-1 \geq r_{1} \geq 0$, $b-2 \geq r_{2} \geq 0$, etc. 
Since $b-b=0$, the largest possible index of the remainder sequence is $b$, so the algorithm terminates in no more than $b$ steps. 
This establishes that the Euclidean Algorithm terminates. 
That each of $G$, $n$, $(r_{1}, \ldots r_{n})$, and $(q_{1},\ldots,q_{n})$ record what is claimed about them in part {\sl (0)} of the theorem statement follows from the description of the terminal step of the algorithm. 

For {\sl (1)}, we induct on $i$ to see that each remainder $r_{i}$ is thusly expressible as an integer linear combination of $a$ and $b$.  
When $i=-1$ we have $r_{-1} = a = (-1)^{0}\cdot{1}\cdot{a} + (-1)^{-1}\cdot{0}\cdot{b}$, and when $i=0$ we have $r_{0} = b = (-1)^{1}\cdot{0}\cdot{a} + (-1)^{0}\cdot{1}\cdot{b}$.   
Assume now that for some $j \in \{0,1,\ldots,n-1\}$, we have $r_{i} = (-1)^{i+1}\overrightarrow{g_{i}}\, a + (-1)^{i}\overrightarrow{h_{i}}\, b$ for each $i \in \{-1,0,\ldots,j\}$. The fact that $r_{j-1} = r_{j}q_{j+1}+r_{j+1}$ implies that 
\begin{eqnarray*}
r_{j+1} & = & r_{j-1} - q_{j+1}r_{j}\\
 & = & \left(\rule[-2mm]{0mm}{4.5mm}(-1)^{j}\overrightarrow{g_{j-1}}\, a + (-1)^{j-1}\overrightarrow{h_{j-1}}\, b\right) + q_{j+1}\left(\rule[-2mm]{0mm}{4.5mm}(-1)^{j}\overrightarrow{g_{j}}\, a + (-1)^{j-1}\overrightarrow{h_{j}}\, b\right)\\
 & = & (-1)^{j}\left(\rule[-2mm]{0mm}{4.5mm}\overrightarrow{g_{j-1}}+\overrightarrow{g_{j}}q_{j+1}\right)a + (-1)^{j-1}\left(\rule[-2mm]{0mm}{4.5mm}\overrightarrow{h_{j-1}}+\overrightarrow{h_{j}}q_{j+1}\right)b\\
 & = & (-1)^{j+2}\overrightarrow{g_{j+1}}\, a + (-1)^{j+1}\overrightarrow{h_{j+1}}\, b,
\end{eqnarray*}
as desired (where the latter equality follows by applying \RecurrenceProp). 

To prove the last claims of part {\sl (1)}, we demonstrate via induction that for each $i \in \{-1,0,\ldots,n+1\}$, we have $r_{n-i} = \overleftarrow{g_{i}}\, G$ and $r_{n-1-i} = \overleftarrow{h_{i}}\, G$, where $G = r_{n-1} = \gcd(a,b)$.  
Well, when $i=-1$ we have $r_{n+1} = G = \overleftarrow{g_{-1}}\, G$ and $r_{n} = 0 = \overleftarrow{h_{-1}}\, G$. 
Now assume that for some $j \in \{-1,0,\ldots,n\}$ it is the case that $r_{n-i} = \overleftarrow{g_{i}}\, G$ and $r_{n-1-i} = \overleftarrow{h_{i}}\, G$ for each $i \in \{-1,0,\ldots,j\}$. 
By this inductive hypothesis, $r_{n-(j+1)} = r_{n-1-j} = h_{j}(q_{n+1-j},\ldots,q_{n})\, G = g_{j+1}(q_{n-j},\ldots,q_{n})\, G$, where the latter equality follows from \RecurrenceProp. 
Similarly, we obtain $r_{n-1-(j+1)} =  r_{n-(j+1)}q_{n-j} + r_{n-j} = q_{n-j}h_{j}(q_{n+1-j},\ldots,q_{n})\, G + g_{j}(q_{n+1-j},\ldots,q_{n})\, G = h_{j+1}(q_{n-j},\ldots,q_{n})\, G$. This completes this induction argument, and the proof of {\sl (1)}.

For part {\sl (2)}, note that if $c|r_{i}$ and $c|r_{i-1}$, then $c|(r_{i-1}q_{i} + r_{i})$, so $c|r_{i-2}$. 
Repeat this reasoning to see that $c|r_{i-3}$, $c|r_{i-4}$, etc, so that eventually we have $c|b$ and $c|a$. 
Since $c$ is a divisor of both $a$ and $b$, then, by part {\sl (1)}, $c\ \rule[-2mm]{0.2mm}{6mm}\left(\rule[-2.5mm]{0mm}{6mm}(-1)^{i+1}\overrightarrow{g_{i}}\, a + (-1)^{i}\overrightarrow{h_{i}}\, b\right)$, so $c$ is a divisor of each $r_{i}$. 
For part {\sl (3)}, since $r_{n}=0$ and $r_{n-1}=G$, then $G|r_{n}$ and $G|r_{n-1}$, and so, by {\sl (2)}, $G|r_{i}$ for each $i \in \{-1,0,\ldots,n\}$. 
In particular, $G$ is a common divisor of $a=r_{-1}$ and $b=r_{0}$, hence, by {\sl (2)}, $\gcd(a,b) \geq G$. 
Since $\gcd(a,b)$ is also a common divisor of $a=r_{-1}$ and $b=r_{0}$, then, again by {\sl (2)}, $\gcd(a,b) | r_{n-1}$, and in paricular $\gcd(a,b) \leq G$. 
So, $\gcd(a,b) = G$. 
Then, every element of the set $\{ax+by\, |\, x,y \in \mathbb{Z}\}$ is a multiple of $G$; this set will consist of all the multiples of $G$ if we can show that $G$ is an integer linear combination of $a$ and $b$. 
Of course, this latter fact is just the $i=n-1$ case of {\sl (1)}.\hfill\QED

{\bf \S \ClosingSection. Closing thoughts.} 
To conclude our discussion, we recall a question from our earlier \EAExample: 
How did we find integers $a$ and $b$ such that the Euclidean Algorithm would terminate in exactly $n=7$ steps, and why $n=7$? 
Well, in seeking examples to illustrate the algorithm to students in Discrete Mathematics courses, we thought that having students perform calculations that involve $g_{7}$ and $h_{7}$ -- which have 13 and 21 terms respectively -- would be sufficient to strongly suggest our proposed connection between the Euclidean Algorithm, B\'{e}zout's Lemma, and our Fibonacci array. 
According to \EATheorem.1, we can find such $a$ and $b$ simply by choosing a sequence $q_{1},q_{2},\ldots,q_{7}$ of quotients and then taking $b:=g_{7}(q_{1},q_{2},\ldots,q_{7})$ and $a:=h_{7}(q_{1},q_{2},\ldots,q_{7})$. 
We leave it as an exercise for the reader to verify that for the sequence $4, 1, 3, 7, 3, 2, 4$ of quotients from \EAExample, we indeed get $b = 935 = g_{7}(4,1,3,7,3,2,4)$ and $a = 4449 = h_{7}(4,1,3,7,3,2,4)$. 
We invite the reader to consider: What aspect of \EATheorem\ allows us to conclude that these numbers $a$ and $b$ are relatively prime?

In \S 3.1 of \cite{Bullynck}, Bullynck observes that certain number theory discourses written separately by Leonhard Euler and Abraham Gotthelf K\"{a}stner in the later 18th century illustrated the Euclidean Algorithm with input integers requiring at most five or six applications of the Division Algorithm in order to return the g.c.d.,  
although K\"{a}stner mentions a case ``where one has to divide 54 times, before one finds the greatest common measure.''  
Bullynck further notes that this latter question was a classic recreational problem originally due to Leonardo Pisano (whose alias we leave it to the reader to discern). 
So, we challenge the reader to use our work here to find the smallest input pair of positive integers $a$ and $b$ (with `smallness' measured, say, by their sum $a+b$) for which the Euclidean Algorithm requires exactly 54 steps to terminate.   
In this case, what is the g.c.d.?

Our effort to better understand the historical record concerning Terquem's problem led us to connect his ideas with Fibonacci numbers, continuants, B\'{e}zout's Lemma, and the Euclidean Algorithm. 
In particular, we have a newfound, or at least now-unforgotten, interpretation of the numbers of the Fibonacci array.
This directly yields combinatorially explicit descriptions of certain polynomials related to continuants. 
And, given positive integers $a$ and $b$, it also directly yields explicit formulas for B\'{e}zout coefficients and for the integers $a/\!\gcd(a,b)$ and $b/\!\gcd(a,b)$. 
All these ideas have antecedents in the literature, but, taken together as presented here, we believe they provide an explicit and combinatorially compelling account of these elementary number theoretic notions and might be advantageously accommodated in introductory discourses on these topics.

%
\renewcommand{\refname}{\large \bf References}
\renewcommand{\baselinestretch}{1.1}

\end{document}